\documentclass[11pt, a4paper]{article}

\usepackage[a4paper, top=2.5cm, bottom=2.5cm, left=2cm, right=2cm]{geometry}
\usepackage[utf8]{inputenc}
\usepackage[T1]{fontenc}
\usepackage[english]{babel}

\usepackage{amsmath,amssymb,amsthm}
\usepackage{mathtools}
\usepackage{booktabs}
\usepackage{enumitem}
\usepackage{hyperref}
\usepackage{comment}

\hypersetup{
    colorlinks=true,
    linkcolor=blue,
    citecolor=blue,
    urlcolor=blue
}

\newtheorem{theorem}{Theorem}[section]
\newtheorem{lemma}[theorem]{Lemma}

\theoremstyle{definition}
\newtheorem{definition}[theorem]{Definition}

\theoremstyle{remark}
\newtheorem{remark}[theorem]{Remark}

\DeclareMathOperator{\dist}{dist}
\DeclareMathOperator{\clm}{clm}
\DeclareMathOperator{\gph}{gph}
\DeclareMathOperator{\conv}{conv}
\DeclareMathOperator{\intt}{int}
\DeclareMathOperator{\bd}{bd}

\title{\textbf{Calmness Modulus of Parameterized Convex Constraint Systems: Exact Closed-Form Formulas via Second-Order Admissibility}}

\author{
  \textbf{J. Camacho}\thanks{Center of Operations Research, Miguel Hern\'{a}ndez University of Elche, 03202 Elche (Alicante), Spain (\texttt{j.camacho@umh.es}). Partially supported by the Spanish Ministry of Science and Innovation (MICINN) grant PID2025-172680NB-I00 and ERDF. This research was also partially supported by the CIGE/2026/47 grant, funded by the Conselleria de Educación, Cultura, Universidades y Empleo, Generalitat Valenciana.}
}
\date{}

\begin{document}

\maketitle

\begin{abstract}
\noindent We resolve the open problem of deriving an exact, closed-form algebraic formula for the calmness modulus of $\mathcal{C}^2$ convex constrained feasible set mappings $\mathcal{F}(b) := \{ x \in \mathbb{R}^n \colon f_t(x) \le b_t, \, \forall t \in T \}$ under right-hand side perturbations at a given pair of the graph $(\bar{b}, \bar{x})$. 

\noindent First, under the condition $0_n \in \intt\left(\conv\left\{ \nabla f_t(\bar{x}) ,\; t \in T(\bar{x}) \right\}\right)$, we prove that the calmness modulus coincides with the inverse distance from the origin to the boundary of the nominal subdifferential. 

\noindent Second, for general $\mathcal{C}^2$ convex systems where the previously mentioned interior condition does not hold, we introduce a second-order admissibility filtering over the active indexes. We show that this filtered family, $\mathcal{D}^0_2(\bar{x})$, yields an exact characterization of the calmness modulus relying solely on nominal data.
\end{abstract}

\noindent \textbf{Keywords:} Calmness modulus, Directional derivatives, Convex analysis.  
\medskip

\noindent \textbf{Mathematics Subject Classification (2020):} 49K40, 90C31, 90C25.

\section{Introduction}

Parametric optimization and variational analysis have proven to be a valuable tool for practitioners that need their solutions to be stable (for a deep introduction in the matter, see e.g. \cite{RockafellarWets98, Mordukhovich06, KlatteKummer02}). Additionally, quantitative stability analysis moves a step further, not only answering to whether a problem/solution enjoys some desired property, but providing concrete measurements that are actionable by decision makers. From an applied mathematical point of view, these quantifiers are essential for other related fields. For instance, they draw a threshold for algorithm convergence (see e.g. \cite{Leventhal08, Lu25, Van25}). Among this rates, the (sharp) Hoffman constant, introduced in the seminal work \cite{Hoffman52} for linear inequality systems, stands up. 

The topic of local quantifiers is quite active, as reflected by the recent literature \cite{Andreani26, Fischer25, Le26}, among others. In particular, this work concerns the so-called calmness modulus in the convex system setup, for which the main antecedents are:
\begin{itemize}
    \item M. J. Cánovas, M. A. López, J. Parra, and F. J. Toledo, \cite{Canovas14}: Characterizing it in the feasible linear inequality system context.
    \item M. J. Cánovas, R. Henrion, M. A. López, and J. Parra, \cite{Canovas16}: Extends the formulae to the $\mathcal{C}^1$ and convex setting, providing upper and lower estimates.
    \item A. Eberhard, V. Roshchina, and T. Sang \cite{Eberhard17}: Extends some of the previous findings to min-max functions. The gap between upper an lower bounds remains.
\end{itemize}

At this point, is worth mentioning that the calmness modulus present a double interest, not only by itself, but also in relation to other measures. Since \cite{Camacho22} proved that  the maximum of certain calmness moduli equals the sharp Hoffman constant in the convex framework, it would result specially useful to obtain an algebraic representation of it.

The current paper closes the gap presented in the literature for an exact characterization of the calmness modulus of of the feasible set mapping of $\mathcal{C}^2$ convex inequality systems under right-hand side perturbations. 

The structure of the paper goes as follows: Section \ref{sec: preliminaries} present the background of the topic, with the needed definitions and state-of-the-art estimates for desired modulus. Section \ref{sec: main} present the main results, Theorems \ref{thm: interior} and \ref{thm: exterior}. The former allows the computation of the modulus under the assumption $0_n \in \intt\left(\conv\left\{ \nabla f_t(\bar{x}) ,\; t \in T(\bar{x}) \right\}\right)$ in the $\mathcal{C}^1$ convex constrained space. In case the assumption falls, the latter result characterizes the quantifier through a second order filtering introduced in Definition \ref{def: second order filtering}. Conclusions and possible future lines of research are gathered in Section \ref{sec: conclusion}. 

\section{Preliminaries}\label{sec: preliminaries}

Let us start by fixing the notation and recalling the main definitions of the paper.
Given $X\subset \mathbb{R}^{p}$, $p\in \mathbb{N}$,
we denote by $\intt\left(X\right)$, $\bd\left(X\right)$, and $\conv\left(X\right)$ the interior,
the boundary, and the convex hull of $X$ respectively.
Hereafter, the convention $\conv\left(\emptyset\right) =\emptyset $ is maintained.

Given a finite index set $T = \{1, \dots, m\}$ and continuous functions $f_t\colon \mathbb{R}^n \to \mathbb{R}$ ($t \in T$), we consider the parameterized feasible set mapping $\mathcal{F}: \mathbb{R}^m \rightrightarrows \mathbb{R}^n$ defined by right-hand side (RHS) perturbations:
\begin{equation*}
\mathcal{F}(b) = \left\{ x \in \mathbb{R}^n \colon f_t(x) \le b_t, \; \forall t \in T \right\}. \label{eq: feasible map}
\end{equation*}

Concerning variational properties, we focus on the \emph{calmness} one, which is considered to be of local nature, since it involves pairs of the graph around a nominal one. The calmness property holds for $\mathcal{F}$ --we particularize the definition for the required multifunction $\mathcal{F}$, but it can be replaced by any multifunction between metric spaces-- at $\left( \bar{b},
\bar{x}\right) \in \gph\mathcal{F}$ (the graph of $\mathcal{F}$)
if there exist a constant $\kappa \geq 0$ and a neighborhood of $\left( 
\bar{b},\bar{x}\right)$, $V\times U$, such that 
\begin{equation}
\dist\left(x,\mathcal{F}\left(\bar{b}\right)\right)\leq \kappa \dist\left( b,\bar{b}\right) 
\text{ for all }x\in \mathcal{F}(b)\cap U\text{ and all }b\in V.
\label{eq: calmness property}
\end{equation}
Along the paper, the \emph{space of variables}, $\mathbb{R}^n$, is equipped with any norm $\Vert\cdot\Vert$, while its dual is defined as standard, $\Vert u\Vert_*\coloneq \max_{\Vert x\Vert=1}\vert u^\top x\vert$. The \emph{parameter space}, $\mathbb{R}^m$, is endowed with the supermum (maximun in our framework, indeed) norm $\Vert \cdot \Vert_\infty$. Their associated distances will be denoted by $\dist(\cdot)$, $\dist_*(\cdot)$ and $\dist_\infty(\cdot)$, respectively.

For completeness, let us recall that such a property can be equivalently written in terms of the inverse function, $\mathcal{F}^{-1}$, i.e.
$b\in \mathcal{F}^{-1}(x)\iff x\in \mathcal{F}(b)$,
appealing to the well known \emph{metric subregularity} property (cf. \cite
[Theorem 3H.3 and Exercise 3H.4]{Dontchev09}) of $\mathcal{F}^{-1}$ at $\left( 
\bar{x},\bar{b}\right)$. Namely, the existence of $\kappa
\geq 0$ and a (possibly smaller) neighborhood $U$ of $\bar{x}$ such
that 
\begin{equation}
\dist\left(x,\mathcal{F}\left(\bar{b}\right)\right)\leq \kappa \dist\left( \bar{b},\mathcal{F}%
^{-1}(x)\right) \text{ for all }x\in U.  \label{eq: met subreg}
\end{equation}

The linear rates $\kappa$ appearing in the variational inequalities \eqref{eq: calmness property} and \eqref{eq: met subreg} are not necessarily uniquely determined. The infimum plausible $\kappa$, known as the \emph{modulus} of the property, is of special interest. It is well known that such moduli coincide for both calmness and metric subregularity, even though their associated neighborhoods might differ.

As a straightforward reformulation of its definition, the calmness modulus can be equivalently written as
\begin{equation}
\clm \mathcal{F}\left(\bar{b}, \bar{x}\right) := \limsup_{\left(b,x\right)\overset{\gph\mathcal{F}}{\longrightarrow }\left( \bar{b},\bar{x}\right) } \frac{\dist\left(x, \mathcal{F}\left(\bar{b}\right)\right)}{\|b-\bar{b}\|_\infty}, \label{eq: calmness mod}
\end{equation}
where $\dist(x, C) := \inf_{z \in C} \|x - z\|$ under the convention $\inf \emptyset\coloneq +\infty$, so that $\dist(x,\emptyset)=+\infty$. Additionally, the agreement $\frac{0}{0}\coloneq 0$ also stands.

Whenever restricted to linear constraints ($f_t(x) = a_t^\top x$), \cite[Theorem 4]{Canovas14} established an exact point-based formula for the calmness modulus. They define the index family $\mathcal{D}(\bar{x})$ as the collection of all subsets $D \subseteq T(\bar{x})\coloneq\{t \in T \colon f_t(\bar{x}) = \bar{b}_t\}$ for which the system in $d \in \mathbb{R}^n$:
\begin{equation}
\left\{ \begin{array}{ll} \nabla f_t(\bar{x})^\top d = 1, & t \in D, \\ \nabla f_t(\bar{x})^\top d < 1, & t \in T(\bar{x}) \setminus D \end{array} \right\} \label{eq: systemD}
\end{equation}
is consistent. Then, the calmness modulus satisfies:
\begin{equation}
\clm \mathcal{F}(\bar{b}, \bar{x}) = \max_{D \in \mathcal{D}(\bar{x})} \left[ \dist_*\left(0_n, \conv\{a_t ,\; t \in D\}\right) \right]^{-1}. \label{eq: clm linear}
\end{equation}

\begin{remark}\label{rem: MFCQ}
Recall that for $\mathcal{C}^1$ inequality systems, the classical \emph{Mangasarian--Fromovitz Constraint Qualification} (MFCQ) holds at $\bar{x}$ if there exists a direction $d \in \mathbb{R}^n$ such that $\nabla f_t(\bar{x})^\top d < 0$ for all $t \in T(\bar{x})$. By Motzkin's Transposition Theorem / Farkas' Lemma, MFCQ is geometrically equivalent to:
\begin{equation*}
0_n \notin \conv\left\{ \nabla f_t(\bar{x}) \;\middle|\; t \in T(\bar{x}) \right\} = \partial f(\bar{x}). \label{eq: MFCQ}
\end{equation*}
Additionally, in the convex setting, MFCQ coincides with the \emph{Slater condition}. When MFCQ holds, $\dist_*(0_n, \partial f(\bar{x})) > 0$, guaranteeing metric regularity (and hence metric subregularity) of the constraint system. 
\end{remark}

In order to extend these findings to $\mathcal{C}^1$ convex systems, we first define the supremum function $f(x) := \max_{t \in T} (f_t(x) - \bar{b}_t)$ and express the calmness modulus via the outer limit of subdifferentials: 
\begin{equation}
\clm \mathcal{F}\left(\bar{b}, \bar{x}\right) = \left[ \dist_*\left(0_n, \limsup_{x \to \bar{x}, f(x) > 0} \partial f(x)\right) \right]^{-1}, \label{eq: clm subdiff}
\end{equation}
In paper \cite{Canovas16}, the authors derive upper and lower bounds to $\limsup_{x \to \bar{x}, f(x) > 0} \partial f(x)$ (and accordingly, to the inverted distance), even though the lower estimate was improved later in \cite{Eberhard17} and extended to min-max functions:

\begin{theorem}\label{thm: sota bounds}
    Given $(\bar{b}, \bar{x}) \in \gph \mathcal{F}$:
    \begin{enumerate}[label=(\alph*)]
        \item[\cite{Canovas16}] We have
        \begin{align}
            \clm \mathcal{F}\left(\bar{b}, \bar{x}\right) &\leq \max_{D\in\mathcal{D}(\bar{x})\cup\mathcal{D}^0(\bar{x})} \left[\dist_*\left(0_n,\conv\left\{\nabla f_t(\bar{x}),\; t\in D\right\}\right)\right]^{-1}\\
            &\leq \left[ \dist_*\left( 0_n, \bd\left(\partial f(\bar{x}) \right)\right) \right]^{-1};
            \label{eq: clm C1 upper Canovas16}
        \end{align}
        where $\mathcal{D}^0$ is the collection of indices that result from interchanging $1$ by $0$ in \eqref{eq: systemD}.
        \item[\cite{Eberhard17}]  If all $f_t$ are convex, we additionally have
        \begin{equation}
            \clm \mathcal{F}\left(\bar{b}, \bar{x}\right) \geq \max_{D\in\mathcal{D}(\bar{x})} \left[\dist_*\left(0_n,\conv\left\{\nabla f_t(\bar{x}),\; t\in D\right\}\right)\right]^{-1},
            \label{eq: clm C1 lower Eberhard17}
        \end{equation}
    \end{enumerate}
\end{theorem}

\section{Main results}\label{sec: main}

In this section, we provide the two main characterizations of the calmness modulus in the context of $\mathcal{C}^1$ and $\mathcal{C}^2$ convex inequality systems under right-hand side perturbations, exposed in Theorems \ref{thm: interior} and \ref{thm: exterior} respectively. To do so, we first observe that the relative position of the origin with respect to the subdifferential of the supremum function plays a fundamental role.

To analyze sequences $x_k \to \bar{x}$ with $f(x_k) > 0$, we decompose vectors into directional components $d \in \mathbb{S}^{n-1} := \{ d \in \mathbb{R}^n \colon \|d\| = 1 \}$. Any sequence $x_k \to \bar{x}$ can be expressed as $x_k = \bar{x} + t_k d_k$, with $t_k \downarrow 0$ and $d_k \to d$. The positive directional cone $S \subseteq \mathbb{S}^{n-1}$ is the set of directions $d \in \mathbb{S}^{n-1}$ for which there exist $t_k \downarrow 0$ and $d_k \to d$ with $f(\bar{x} + t_k d_k) > 0$ for all $k \in \mathbb{N}$. Define the first-order directional value function and active directional indices:
\begin{equation}
M(d) := \max_{t \in T(\bar{x})} \nabla f_t(\bar{x})^\top d, \qquad D(d) := \left\{ t \in T(\bar{x}) \colon \nabla f_t(\bar{x})^\top d = M(d) \right\}. \label{eq: M and D}
\end{equation}

The unit sphere $\mathbb{S}^{n-1}$ partitions into three disjoint directional sets based on $M(d)$:
\begin{itemize}
  \item Strictly Positive Directional Set: $S^+ := \left\{ d \in \mathbb{S}^{n-1} \colon M(d) > 0 \right\}$.
  \item Critical Directional Set: $S^0 := \left\{ d \in \mathbb{S}^{n-1} \colon M(d) = 0 \right\}$.
  \item Strictly Negative Directional Set: $S^- := \left\{ d \in \mathbb{S}^{n-1} \colon M(d) < 0 \right\}$.
\end{itemize}

\begin{lemma}\label{lem: directional cones}
The directional cones satisfy $S^+ \subseteq S$, $S^- \cap S = \emptyset$, and $S = S^+ \cup (S^0 \cap S)$.
\end{lemma}

\begin{proof}
The first two are straightforward from the definitions. Aiming to prove the last equality, let us consider the first-order Taylor expansion $f_t(\bar{x} + \varepsilon d) = \bar{b}_t + \varepsilon \nabla f_t(\bar{x})^\top d + o(\varepsilon)$. Taking the maximum over $t \in T(\bar{x})$ gives $f(\bar{x} + \varepsilon d) = \varepsilon M(d) + o(\varepsilon)$. 

\smallskip
\noindent Case $d \in S^+$: We have $M(d) > 0$ and hence, for sufficiently small $\varepsilon > 0$, $f(\bar{x} + \varepsilon d)>0$. Thus, $d \in S$.

\smallskip
\noindent 
Case $d \in S^-$: Analogously, we have $M(d) < 0$, which yields   $f(\bar{x} + \varepsilon d)<0$ and so $d \notin S$.

Combining these yields $S = S^+ \cup (S^0 \cap S)$.
\end{proof}

\begin{remark}\label{rem: 1-to-1}
    Note that there is a 1-to-1 link between $S^+$ and $\mathcal{D}(\bar{x})$, as well as, $S^0$ with $\mathcal{D}^0(\bar{x})$. 
    
    If $d \in S^+$ ($M(d) > 0$), scaling $d' = \frac{d}{M(d)}$ yields $\nabla f_t(\bar{x})^\top d' = 1$ for $t \in D(d)$ and $\nabla f_t(\bar{x})^\top d' < 1$ for $t \in T(\bar{x}) \setminus D(d)$. Thus, $D(d) \in \mathcal{D}(\bar{x})$.
    
    If $d \in S^0$ ($M(d) = 0$), then $\nabla f_t(\bar{x})^\top d = 0$ for $t \in D(d)$ and $\nabla f_t(\bar{x})^\top d < 0$ for $t \in T(\bar{x}) \setminus D(d)$. Thus, $D(d) \in \mathcal{D}^0(\bar{x})$.
    \end{remark}

As previously announced, the position of the origin with respect to the subdifferential induces different behaviors in $S^+$, $S^0$ and $S^-$. This fact motivates the separated studies.

\begin{theorem}\label{thm: interior}
Let $(\bar{b}, \bar{x}) \in \gph \mathcal{F}$, and assume all constraint functions $f_t$ ($t \in T$) are convex and $\mathcal{C}^1$-smooth. If $0_n \in \intt\left(\conv\left\{ \nabla f_t(\bar{x}) ,\; t \in T(\bar{x}) \right\}\right)$,
then $S^0 = S^-= \emptyset$, $S^+ = \mathbb{S}^{n-1}$, and the calmness modulus is given by:
\begin{equation*}
\clm \mathcal{F}\left(\bar{b}, \bar{x}\right) = \max_{D\in\mathcal{D}(\bar{x})} \left[\dist_*\left(0_n,\conv\left\{\nabla f_t(\bar{x}),\; t\in D\right\}\right)\right]^{-1} = \left[ \dist_*\left( 0_n, \bd\left(\partial f(\bar{x}) \right)\right) \right]^{-1}. \label{eq: clm C1}
\end{equation*}
\end{theorem}

\begin{proof}
Since $f_t$ are $\mathcal{C}^1$ convex for all $t\in T$ , $f(x) = \max_{t \in T}(f_t(x) - \bar{b}_t)$ is also convex and continuous. Its subdifferential at $\bar{x}$ is the compact convex set $\partial f(\bar{x}) = \conv\left\{ \nabla f_t(\bar{x}) ,\, t \in T(\bar{x}) \right\}$. The support function of $\partial f(\bar{x})$ for any unit direction $d \in \mathbb{S}^{n-1}$ is precisely $ M(d) = \max_{t \in T(\bar{x})} \nabla f_t(\bar{x})^\top d $.

By assumption, $0_n \in \intt(\partial f(\bar{x}))$, which implies that the support function is strictly positive in every direction \cite[Theorem 13.1]{Rockafellar70}:
$M(d) > 0$ for all $d \in \mathbb{S}^{n-1}$. Thus, $S^+ = \mathbb{S}^{n-1}$ and $S^0 = S^-= \emptyset$.

Now take any direction $d \in \mathbb{S}^{n-1}$ and sequence $t_k \downarrow 0$. For all small $\varepsilon > 0$, non-active constraints $s \notin T(\bar{x})$ satisfy $f_s(\bar{x} + \varepsilon d) - \bar{b}_s < 0$. For active constraints $t \in T(\bar{x})$, $f_t(\bar{x} + \varepsilon d) - \bar{b}_t = \varepsilon \nabla f_t(\bar{x})^\top d + o(\varepsilon)$. Thus, $T(\bar{x} + \varepsilon d) = D(d)$ for all sufficiently small $\varepsilon > 0$.
Since $f_t \in \mathcal{C}^1$, continuity of gradients yields:
\begin{equation*}
    \lim_{\varepsilon \downarrow 0} \partial f(\bar{x} + \varepsilon d) = \conv\left\{ \nabla f_t(\bar{x}) ,\, t \in D \right\},
\end{equation*}
Taking the outer limit over all sequences $x_k \to \bar{x}$ with $f(x_k) > 0$ (which correspond to rays $d \in S^+$):
\begin{equation}
\limsup_{x \to \bar{x}, f(x) > 0} \partial f(x) = \bigcup_{d \in S^+} \conv\left\{ \nabla f_t(\bar{x}) ,\,  t \in D(d) \right\}=\bd\left(\partial f(\bar{x})\right). \label{eq: outer limit subdiff 1}
\end{equation}
Taking distances to the origin in \eqref{eq: outer limit subdiff 1}, together with formula \eqref{eq: clm subdiff} and Remark \ref{rem: 1-to-1} finishes the proof.
\end{proof}

When $0_n \notin \intt(\partial f(\bar{x}))$, the critical set of directions $S^0$ (or equivalently $\mathcal{D}^0(\bar{x})$) may not by empty, so the previous argument cannot be replicated. According to Theorem \ref{thm: sota bounds}, the calmness modulus is bounded from below and above via $\mathcal{D}(\bar{x})$ and $\mathcal{D}(\bar{x})\cup \mathcal{D}^0(\bar{x})$ respectively. Moreover, depending on the system, it can be the case that either extremes or anything in between is the exact value, as exposed in the subsection \ref{sub: example} and also commented in \cite[Example 3.2, Example 3.3, Remark 3.6]{Canovas16}.

\begin{definition}\label{def: second order filtering}
Let $f_t$ ($t \in T$) be $\mathcal{C}^2$-smooth convex functions. 
The family of \emph{second-order index sets}, denoted by $ \mathcal{D}^0_2(\bar{x})$, is the collection of all sets $D\in \mathcal{D}^0(\bar{x})$ for which there exists a direction $d \in \mathbb{S}^{n-1}$ satisfying:
\begin{enumerate}[label=\arabic*.]
  \item $\nabla f_t(\bar{x})^\top d = 0$ for all $t \in D$, and $\nabla f_s(\bar{x})^\top d < 0$ for all $s \in T(\bar{x}) \setminus D$.
  \item $D = \left\{ t \in D(d) \;\middle|\; d^\top \nabla^2 f_t(\bar{x}) d = \max_{s \in D(d)} d^\top \nabla^2 f_s(\bar{x}) d > 0 \right\}$
\end{enumerate}
\end{definition}

\begin{theorem}\label{thm: exterior}
Let $(\bar{b}, \bar{x}) \in \gph \mathcal{F}$. Assume, for all $t \in T$, that $f_t$ are convex and $\mathcal{C}^2$-smooth. Then:
\begin{enumerate}[label=(\alph*)]
  \item If $0_n \in \conv\{\nabla f_t(\bar{x}) ,\, t \in D^0\}$ for some $D^0 \in \mathcal{D}^0_2(\bar{x})$, then $\clm \mathcal{F}(\bar{b}, \bar{x}) = +\infty$.
  \item Otherwise, the calmness modulus is given by
  \begin{equation}
  \clm \mathcal{F}\left(\bar{b}, \bar{x}\right) = \max_{D \in \mathcal{D}(\bar{x}) \cup \mathcal{D}^0_2(\bar{x})} \left[ \dist_*\left(0_n, \conv\{\nabla f_t(\bar{x}) ,\, t \in D\}\right) \right]^{-1}<+\infty. \label{eq: clm exact}
  \end{equation}
\end{enumerate}
\end{theorem}

\begin{proof}
By Lemma \ref{lem: directional cones}, $S = S^+ \cup (S^0 \cap S)$. We examine outer limits of subdifferential along sequences $x_k = \bar{x} + \varepsilon_k d_k \to \bar{x}$ with $f(x_k) > 0$ by separating directions into $S^+$ and $S^0 \cap S$:

\smallskip
\noindent Case $d \in S^+$: As shown in Theorem \ref{thm: interior}, $M(d) > 0$ ensures $f(\bar{x} + \varepsilon d) > 0$ for all small $\varepsilon > 0$. The subdifferential limit set generated along such rays is $\conv\{\nabla f_t(\bar{x}) ,\, t \in D(d)\}$. By Remark \ref{rem: 1-to-1}, as $d$ ranges over $S^+$, the index sets $D(d)$ cover $\mathcal{D}(\bar{x})$.

\smallskip
\noindent Case $d \in S^0 \cap S$: Since $d \in S^0$, $M(d) = 0$. For non-maximal active constraints $s \in T(\bar{x}) \setminus D(d)$, $\nabla f_s(\bar{x})^\top d < 0$, so $f_s(\bar{x} + \varepsilon d) - \bar{b}_s < 0$ for small $\varepsilon > 0$. For maximal active constraints $t \in D(d)$, we have $f_t(\bar{x}) - \bar{b}_t = 0$ and $\nabla f_t(\bar{x})^\top d = M(d) = 0$. Hence the 0-th and 1-st order Taylor terms vanish, yielding:
\begin{equation*}
f_t(\bar{x} + \varepsilon d) - \bar{b}_t = \frac{1}{2} \varepsilon^2 d^\top \nabla^2 f_t(\bar{x}) d + o(\varepsilon^2), \quad \forall t \in D(d). \label{eq: Taylor second order}
\end{equation*}
Therefore, $f(\bar{x} + \varepsilon d) = \max_{t \in D(d)} (f_t(\bar{x} + \varepsilon d) - \bar{b}_t) > 0$ holds for small $\varepsilon > 0$ if and only if:
\begin{equation*}
    \max_{t \in D(d)} d^\top \nabla^2 f_t(\bar{x}) d > 0.
\end{equation*}

Furthermore, an index $t \in D(d)$ belongs to the active set $T(\bar{x} + \varepsilon d)$ for small $\varepsilon > 0$ if and only if its second-order rate attains this maximal value:
\begin{equation*}
d^\top \nabla^2 f_t(\bar{x}) d = \max_{s \in D(d)} d^\top \nabla^2 f_s(\bar{x}) d.
\end{equation*}
Indices in $D(d)$ with strictly smaller second-order values satisfy $f_s(\bar{x} + \varepsilon d) < f(\bar{x} + \varepsilon d)$ and remain strictly inactive for small $\varepsilon > 0$. Consequently, subdifferentials along rays $d \in S^0 \cap S$ converge exclusively to sets $\conv\{\nabla f_t(\bar{x}) ,\, t \in D^0\}$ corresponding to index sets $D^0 \in \mathcal{D}^0_2(\bar{x})$ as defined in Definition \ref{def: second order filtering}.

\smallskip
\noindent Combining both cases, the outer limit of subdifferentials is:
\begin{equation}
\limsup_{x \to \bar{x}, f(x) > 0} \partial f(x) = \bigcup_{D \in \mathcal{D}(\bar{x}) \cup \mathcal{D}^0_2(\bar{x})} \conv\left\{ \nabla f_t(\bar{x}) ,\, t \in D \right\}. \label{eq: outer limit subdiff 2}
\end{equation}

If $0_n \in \conv\{\nabla f_t(\bar{x}) ,\, t \in D^0\}$ for some $D^0 \in \mathcal{D}^0_2(\bar{x})$, then $\dist_*(0_n, \limsup_{x \to \bar{x}, f(x) > 0} \partial f(x)) = 0$, so \eqref{eq: clm subdiff} yields $\clm \mathcal{F}(\bar{b}, \bar{x}) = +\infty$, establishing part (a).
Otherwise, taking dual distances in \eqref{eq: outer limit subdiff 2} and applying \eqref{eq: clm subdiff} yields \eqref{eq: clm exact}, establishing part (b).
\end{proof}

If $f_t$ for $t\in T$ are affine, the second order family $\mathcal{D}^0_2(\bar{x})$ is empty, so that Theorem \ref{thm: exterior} recovers formula \eqref{eq: clm linear}.

\subsection{Illustrative example}\label{sub: example}

To highlight the significance of the second-order filtering, we present a low-dimensional example where prior bounds from \cite{Canovas16, Eberhard17} present a gap, whereas Theorem~\ref{thm: exterior} yields the exact calmness modulus.

Consider the decision space $\mathbb{R}^2$ equipped with the Euclidean norm $\|\cdot\|_2$, and parameter space $\mathbb{R}^4$ with the Chebyshev norm $\|\cdot\|_\infty$. Let $\bar{x} = (0,0)^\top \in \mathbb{R}^2$ and nominal parameter $\bar{b} = (0,0,0,0)^\top \in \mathbb{R}^4$. We define five convex constraint functions $f_t \colon \mathbb{R}^2 \to \mathbb{R}$ ($t \in T = \{1, 2, 3, 4\}$):
\begin{align*}
f_1(x_1, x_2) &:= x_1 &\implies \nabla f_1(\bar{x}) &= \begin{pmatrix} 1 \\ 0 \end{pmatrix}, & \nabla^2 f_1(\bar{x}) &= \begin{pmatrix} 0 & 0 \\ 0 & 0 \end{pmatrix}, \\ 
f_2(x_1, x_2) &:= \frac{1}{2}x_1 + 2x_2^2               &\implies \nabla f_2(\bar{x}) &= \begin{pmatrix} 1/2 \\ 0 \end{pmatrix},   & \nabla^2 f_2(\bar{x}) &= \begin{pmatrix} 0 & 0 \\ 0 & 4 \end{pmatrix}, \\ 
f_3(x_1, x_2) &:= \frac{1}{3}x_1 + x_2^2               &\implies \nabla f_3(\bar{x}) &= \begin{pmatrix} 1/3 \\ 0 \end{pmatrix},   & \nabla^2 f_3(\bar{x}) &= \begin{pmatrix} 0 & 0 \\ 0 & 2 \end{pmatrix}, \\ 
f_4(x_1, x_2) &:= x_1 + x_2                &\implies \nabla f_4(\bar{x}) &= \begin{pmatrix} 1 \\ 1 \end{pmatrix},  & \nabla^2 f_4(\bar{x}) &= \begin{pmatrix} 0 & 0 \\ 0 & 0 \end{pmatrix}, 
\end{align*}
All constraints are active at $\bar{x} = (0,0)^\top$, so $T(\bar{x}) = \{1, 2, 3, 4\}$. It is not difficult to check that $\{1,4\}\in\mathcal{D}(\bar{x})$ and $\{1,2,3\}\in\mathcal{D}^0(\bar{x})$ are maximal with respect to the inclusion order. Accordingly, formulae \eqref{eq: clm C1 upper Canovas16} and \eqref{eq: clm C1 lower Eberhard17} provide the bounds
\begin{equation*}
    1=\dist_2\left(0_2,\begin{pmatrix} 1 \\ 0 \end{pmatrix}\right)^{-1}\leq \clm \mathcal{F}\left(\bar{b}, \bar{x}\right) \leq \dist_2\left(0_2,\begin{pmatrix} 1/3 \\ 0 \end{pmatrix} \right)^{-1}=3.
\end{equation*}

For the direction $d = (0, -1)^\top \in \mathbb{S}^1$, we have $D(d) = \{1, 2, 3\}$. Evaluating the second order rates
\begin{equation*}
    d^\top \nabla^2 f_1(\bar{x}) d = 0, \quad d^\top \nabla^2 f_2(\bar{x}) d = 4, \quad d^\top \nabla^2 f_3(\bar{x}) d = 2,
\end{equation*}
we obtain $\{2\}\in \mathcal{D}_2^0(\bar{x})$, setting up  
\begin{equation*}
    \clm \mathcal{F}\left(\bar{b}, \bar{x}\right)=\dist_2\left(0_2,\begin{pmatrix} 1/2 \\ 0 \end{pmatrix} \right)^{-1}=2.
\end{equation*}

\section{Conclusion}\label{sec: conclusion}

This paper resolves the open problem of determining the exact calmness modulus for feasible set mappings of $\mathcal{C}^2$ convex inequality systems under right-hand side perturbations. By introducing a second-order admissibility filtering over critical active index sets, we eliminate non-attainable boundary subdifferentials and derive a closed-form formula relying solely on nominal first- and second-order data. Future work includes extending these second-order criteria to non-convex variational settings, as well as relaxing the $\mathcal{C}^2$-smoothness necessity in the convex framework.

\section*{Declarations}

\noindent Data availability is not applicable to this article as no new datasets were created or analyzed in this study.

\noindent The authors declare that they have no conflict of interest of any kind related to the manuscript.

\bibliographystyle{elsarticle-num}

\end{document}